\documentclass[12pt,twoside,a4paper]{article}
\usepackage{article}
\usepackage{notation}

\newcommand{\papertitle}{A Semi-equivariant Dixmier-Douady Invariant}
\newcommand{\paperauthor}{Simon Kitson}

\begin{document}
\color{textcolor}

\pagestyle{fancy}
\fancyhead{}
\fancyhead[EC]{\fontsize{9}{12}\selectfont{\MakeUppercase{\paperauthor}}}
\fancyhead[OC]{\fontsize{9}{12}\selectfont{\MakeUppercase{\papertitle}}}
\fancyhead[EL]{\fontsize{9}{12}\selectfont{\thepage}}
\fancyhead[OR]{\fontsize{9}{12}\selectfont{\thepage}}
\fancyfoot{}
\renewcommand{\headrulewidth}{0pt}
\renewcommand{\footrulewidth}{0pt}

\title{\normalsize\bf\MakeUppercase{\papertitle}}
\author{\fontsize{10}{12}\selectfont\MakeUppercase{\paperauthor}\thanks{I would like to thank the Mathematical Sciences Institute of the Australian National University for the postdoctoral fellowship which supported this research.}}
\date{\vspace{-1cm}}
\maketitle

\begin{center}
    \begin{minipage}{0.875\textwidth}
    {\scshape\noindent Abstract.}{\ \small
%
%
A generalisation of 
    the equivariant Dixmier-Douady invariant 
    is constructed as a 
    second-degree cohomology class
    within a new
    semi-equivariant \v{C}ech cohomology theory.
    This invariant obstructs liftings
    of semi-equivariant principal bundles
    that are 
    associated to central exact sequences 
    of structure groups
    in which each structure group is 
    acted on by the equivariance group.
The results and methods 
    described can be applied
    to the study of complex vector bundles 
    equipped with linear/anti-linear actions,
    such as Atiyah's Real vector bundles.

}
    \end{minipage}
\end{center}
\begin{center}
\begin{minipage}{0.875\textwidth}
\tableofcontents 
\end{minipage}
\end{center}
%
%

\section{Introduction}


A \emph{Real vector bundle} $(E,\tau)$
    is a complex vector bundle equipped 
    with an anti-linear involution
    that covers an involution on its base space
    \cite{Atiyah.1966.K-Theory-and-Reality}.
    The $\Unitary(n)$-frame bundle $\Fr(E)$
    of a Real vector bundle
    is equipped with two actions:
    a left action of
    $\Z_2$
    induced by $\tau$, and 
    a right action of $\Unitary(n)$.
    Due to the anti-linearity of $\tau$, 
    these actions do not commute.
    Rather, they combine into an action
    of $\Z_2 \sdp \Unitary(n)$,
    where $\Z_2$ acts on $\Unitary(n)$
    by elementwise conjugation.

More generally,
    if $G$ is a $\Gam$-group%
        \footnote{%
        a group equipped with 
        an action of $\Gam$ 
        by group automorphisms.}
    and $P$
    is a
    principal $G$-bundle
    equipped with a left action of $\Gam$
    that maps fibres to fibres and 
    satisfies
    $\gam(pg) = (\gam p)(\gam g)$
    for all
        $\gam \in \Gam$, 
        $p \in P$ and
        $g \in G$,
    then
    the actions on $P$
    combine into an action of $\Gam \sdp G$.
    In this situation,
    $P$ is described as a
    \emph{$\Gam$-semi-equivariant 
          principal $G$-bundle}.
    When the $\Gam$-action on $G$ is trivial, 
    $P$ is an equivariant principal bundle
    in the usual sense.

This paper solves the following 
    lifting problem
    for semi-equivariant principal bundles
    (see Theorem \ref{thm:mainthm}):

    \begin{minipage}[c]{0.9\textwidth}
    \emph{%
        Given             
        a central short exact sequence 
            $
                A \os{\al}{\rightarrow}
                B \os{\be}{\rightarrow}
                C
            $
        of $\Gam$-groups
        and 
        a $\Gam$-semi-equivariant 
        principal $C$-bundle $P$,
        classify 
        the liftings of $P$ by $\be$
        to 
        a $\Gam$-semi-equivariant 
        principal $B$-bundle.
    }
    \end{minipage}\\

\noindent
In particular, 
    the obstruction to such liftings
    is identified as a
    \emph{semi-equivariant Dixmier-Douady 
          invariant}.
    This new invariant lies in a new
    \emph{semi-equivariant \v{C}ech cohomology 
          theory},
    which will be constructed 
    in \S \ref{sec:semiCo}.
    The semi-equivariant 
    Dixmier-Douady invariant generalises 
    the equivariant Dixmier-Douady invariant,
    which lies in equivariant cohomology.
The relationship of this method 
    to existing work, and 
    its possible applications
    will be discussed 
    in \S\ref{sec:relwork}.


\comment{
    The usual notion of a
        $\Gam$-equivariant principal $G$-bundle
        consists of
        a principal $G$-bundle $P$
        equipped with an action of $\Gam$
        that maps fibres to fibes
        and 
        commutes with the principal $G$-action.
        The isomorphism classes of  
        such bundles
        correspond to 
        equivalence classes 
        of equivariant transition cocycles
        $\phi \in \TCocycles^1_\Gam(X,G)$.
    Associated to 
        a cocycle $\phi \in \TCocycles^1_\Gam(X,G)$
        and 
        a central exact sequence
        \begin{equation}
        \label{eq:equiexact}
            1 
            \rightarrow 
            A 
            \rightarrow 
            H 
            \os{\vphi}{\rightarrow} 
            G 
            \rightarrow 
            1
        \end{equation}
        is a lifting problem:
            \emph{find all classes
            $\psi \in \TCocycles^1_\Gam(X,H)$}
            \emph{such that} $\vphi(\psi) = \phi$.
        A theorem of Dixmier and Douady 
        resolves this problem
        by constructing an exact sequence
        \begin{align}
        \label{eq:equiHexact}
        H_\Gam^1(X,A)     
          \os{}{\rightarrow}
        \TCocycles_\Gam^1(X,H)
          \os{\vphi_*}{\rightarrow}
        \TCocycles_\Gam^1(X,G)
          \os{\Del}{\rightarrow}
        H_\Gam^2(X,A),
        \end{align}
        where $H_\Gam^p(X,A)$
        denotes the equivariant \v{C}ech cohomology.
        The class $\Del(\phi) \in H_\Gam^2(X,A)$
        is the equivariant 
        \emph{Dixmier-Douady invariant}
        of $\phi$.
        When $\Del(\phi) \neq 0$,
        there are no 
        liftings 
        $\psi \in \TCocycles^1_\Gam(X,H)$
        such that 
        $\vphi(\psi) = \phi$.
        When $\Del(\phi) = 0$,
        the liftings $\psi$ 
        correspond non-canonically 
        to the classes of $H_\Gam^1(X,A)$.    
        For this reason, 
        the Dixmier-Douady invariant
        is often described 
        as the \emph{obstruction} 
        to the existence of a liftings.
    
    By generalising this entire argument
        to the semi-equivariant setting
        it is possible to resolve the 
        semi-equivairnat lifting problem.
}



\section{Semi-equivariant Principal Bundles}
\label{ch:semi}




Before examining 
    semi-equivariant principal bundles, 
    the notion of a semi-direct product 
    is breifly reviewed.

\begin{defn}
Let $\Gam$ be a Lie group.
    A \emph{(smooth) $\Gam$-group} $(G,\tht)$ 
    is a Lie group 
    equipped with a smooth action
    \begin{equation*}
    \tht: \Gam \rightarrow \Aut(G).
    \end{equation*}
A \emph{homomorphism} 
    $\vphi: G \rightarrow H$ of $\Gam$-groups 
    is a homomorphism of Lie groups 
    such that, 
    for $\gam \in \Gam$ and $g \in G$,
    \begin{equation}
        \vphi(\gam g) = \gam\vphi(g).
    \end{equation}
\end{defn}

\begin{defn}
    Let $(G,\tht)$ be a $\Gam$-group.
    The \emph{(outer) semi-direct product} 
        $\Gam \sdp_\tht G$ 
    is the Lie group consisting 
    of elements 
       $(\gam,g) \in \Gam \times G$ 
    with multiplication defined, 
    for $\gam_i \in \Gam$ and $g_i \in G$,
    by
    \begin{equation*}
        (\gam_1,g_1)(\gam_2,g_2) 
            := 
        (\gam_1\gam_2, g_1(\gam_1 g_2)).
    \end{equation*}
\end{defn}

One situation in which
    semi-direct product groups
    arise is when $G$ and $\Gam$ 
    both act on an object $X$
    and satsify the relation 
        $\gam (gx) = (\gam g) (\gam x)$,
    for some action $\tht$ of $\Gam$ on $G$.
    In this case, 
    the two actions combine 
    to form a single action
    of the group 
        $\Gam \sdp_\tht G$ 
    by $(\gam,g)x := g(\gam x)$.
\begin{exam}
\label{exam:UKSDP}
The standard $\Unitary(1)$-action on $\C$
    and the $\Z_2$-action on $\C$ by conjugation,
    combine into a
        $\Z_2 \sdp_\kap \Unitary(1)$-action
    on $\C$,
    where $\kap$ is the $\Z_2$-action 
    on $\Unitary(1)$
    by conjugation.
\end{exam}


Semi-equivariant principal bundles
    generalises
    equivariant principal bundle
    by using a $\Gam$-group $(G,\tht)$
    as the structure group.
    The action $\tht$ determines 
    the commutation relation between the 
    left action of $\Gam$ and right action of $G$
    on the total space of the principal bundle.
    These actions combine into
    an action of the semi-direct product 
        $\Gam \sdp_\tht G$.
In the following definitions, 
    let $(G,\tht)$ be a smooth $\Gam$-group
    and $X$ be a manifold equipped 
    with a smooth $\Gam$-action.
\begin{defn}
A \emph{(smooth) $\Gam$-semi-equivariant 
        principal $(G,\tht)$-bundle}
    over $X$
    is a smooth principal $G$-bundle  
        $\pi: P \rightarrow X$
    equipped with a smooth left action of $\Gam$
    such that, 
    for $\gam \in \Gam$, $p \in P$ and $g \in G$,
    \begin{align*}
            \pi(\gam p) &= \gam \pi(p) &
            \gam(pg) &= (\gam p)(\gam g).
    \end{align*}
\end{defn}

\begin{defn}
An \emph{isomorphism} $\vphi: P \rightarrow Q$
    of 
    $\Gam$-semi-equivariant 
    principal $(G,\tht)$-bundles
    is a diffeomorphism
    such that,
    for $\gam \in \Gam$, $p \in P$ and $g \in G$,
    \begin{align*}
        \pi_P &= \pi_Q \circ \vphi &
        \vphi(pg) &= \vphi(p) g &
        \vphi(\gam p) &= \gam \vphi(p).
    \end{align*}
\end{defn}

Next, let 
    $\lam: (G,\tht) \rightarrow (H,\vtht)$ 
    be a homomorphism of $\Gam$-groups, and 
    $Q$ be 
    a $\Gam$-semi-equivariant 
    principal $(H,\vtht)$-bundle.
\begin{defn}
A \emph{lifting} of $Q$ by $\lam$ 
    is a pair $(P,\vphi)$,
    where 
    $P$ is 
    a $\Gam$-semi-equivariant 
    principal $(G,\tht)$-bundle 
    and 
        $\vphi: P \rightarrow Q$ 
    is a smooth map
    such that, 
    for $\gam \in \Gam$, $p \in P$ and $g \in G$,
    \begin{align*}
        \pi_P &= \pi_Q \circ \vphi &
        \vphi(pg) &= \vphi(p) \lam(g) &
        \vphi(\gam p) &= \gam \vphi(p).
    \end{align*}
\end{defn}

\begin{defn}
Two liftings
        $(P_1, \vphi_1)$ and $(P_2,\vphi_2)$ 
    of $Q$ by $\lam$ are \emph{equivalent} 
    if there is an isomorphism
       $\psi: P_1 \rightarrow P_2$
  such that $\vphi_2 \circ \psi = \vphi_1$.
\end{defn}

The set of smooth $\Gam$-semi-equivariant 
    principal $(G,\tht)$-bundles
    will be denoted $\PBundles_\Gam(X,(G,\tht))$, and the
    isomorphisms classes will be denoted 
        $\PBundles^{\iso}_\Gam(X,(G,\tht))$.

\section{Semi-equivariant Transition Cocycles}
\label{sec:SETC}

Transition cocycles are used 
    to extract global topological information 
    from a principal bundle
    into a form which is more easily analysed.
A transition cocycle 
    over an open cover $\cU := \set{U_a}$ 
    with values in a Lie group $G$
    is a collection of smooth maps 
    $\phi_a: U_a \rightarrow G$.
    Maps on overlapping open sets 
    are required to satisfy 
    a \emph{cocycle condition}.
    This condition ensures that 
    the cocycle can be used 
    to glue together 
    the patches $U_a \times G$ 
    into a principal $G$-bundle.

In the equivariant setting, 
    a transition cocycle consists of maps 
        $\phi_a(\gam, \cdot): U_a \rightarrow G$
    for each $U_a \in \cU$ and $\gam \in \Gam$.
    The equivariant cocycle condition 
    then ensures that 
    the elements $\phi_a(1, \cdot)$ 
    can be used construct the total space 
    of a principal $G$-bundle, 
    and that the elements $\phi_a(\gam,\cdot)$ 
    can be used to construct a $\Gam$-action.
    The derivation 
    of the equivariant cocycle condition 
    uses the fact that 
    the left and right actions commute.

Semi-equivariant transition cocycles 
    can be defined in a similar fashion 
    to equivariant transition cocycles.
    However, the left and right actions 
    on a $\Gam$-semi-equivariant 
    principal $(G,\tht)$-bundle
    form an action of $\Gam \sdp_\tht G$.
    Thus, the commutation relation 
    between the left and right actions 
    is controlled by $\tht$, and 
    the action $\tht$ 
    appears in 
    the semi-equivariant cocycle condition.
    When this cocycle condition is satisfied, 
    the elements        
        $\phi_a(1, \cdot)$ 
    in a cocycle 
    can be used to construct 
    the total space 
    of a semi-equivariant principal bundle, 
    and the elements 
        $\phi_a(\gam,\cdot)$ 
    can be used to construct 
    a semi-equivariant $\Gam$-action.

Throughout this section, let 
    $X$ be a $\Gam$-space,
    $(G,\tht)$ be a $\Gam$-group and 
    $\cU := \set{U_a}$ be an open cover of $X$.
    The cover $\cU$ is not required to be invariant.

\begin{defn}
\label{defn:SECocycle}
A \emph{(smooth) $\Gam$-semi-equivariant 
        $(G,\tht)$-valued transition cocycle}
    over $\cU$
    is a collection of smooth maps
    \begin{equation*}
          \phi 
            := 
          \set{
              \phi_{ba}(\gam,\cdot) 
                : 
              U_a \cap \gam^\inv U_b \rightarrow G 
                \mid
              U_a \cap \gam^\inv U_b \neq \emptyset
              }, 
    \end{equation*}
    satisfying
    \begin{align}\label{eq:TransitionCocycleProperty}
       \phi_{aa}(1,x_0) &= 1
        &
       \phi_{ca}(\gam'\gam,x) 
            &= 
        \phi_{cb}(\gam',\gam x)(\gam' \phi_{ba}(\gam,x)),
    \end{align}
    for 
    $x_0 \in U_a$,
    $\gam',\gam \in \Gam$
    and
    $
        x 
        \in 
        U_a 
            \cap 
        \gam^\inv U_b 
            \cap 
        (\gam'\gam)^\inv U_c
    $.
\end{defn}

Note that the conditions
        \eqref{eq:TransitionCocycleProperty} 
    define a non-equivariant cocycle
    when restricted to $\gam = 1$, and
    an equivariant cocycle when $\tht = \id$.

\begin{defn}
\label{defn:SETCEquiv}
    An \emph{equivalence} of
    $\Gam$-semi-equivariant $(G,\tht)$-valued 
    transition cocycles
    $\phi^1$ and $\phi^2$ 
    with cover $\cU$ 
    is a collection of smooth maps 
    \begin{equation*}
        \mu := \set{\mu_a : U_a \rightarrow G}
    \end{equation*}
    such that
    \begin{equation*}
        \mu_b(\gam x) \phi^1_{ba}(\gam,x) 
            = 
        \phi^2_{ba}(\gam,x) (\gam \mu_a(x)),
    \end{equation*}
    for 
        $\gam \in \Gam$ and 
        $x \in U_a \cap \gam^\inv U_b$ .
\end{defn}

Next, let 
        $\lam: (G,\tht) \rightarrow (H,\vtht)$ 
    be a homomorphism of $\Gam$-groups, and 
    $\phi$ be a $\Gam$-semi-equivariant 
    $(H,\vtht)$-valued transition cocycle 
    over $\cU$.

\begin{defn}
A \emph{lifting} of $\phi$ by $\lam$ 
    is a $\Gam$-semi-equivariant 
    $(G,\tht)$-valued transition cocycle $\psi$
  such that
    $\lam \circ \psi_{ba} = \phi_{ba}$.
\end{defn}

\begin{defn}
Two liftings $\psi^1$ and $\psi^2$ of $\phi$ by $\lam$
    are \emph{equivalent} if there exists 
    an equivalence $\mu$
    between $\psi^1$ and $\psi^2$. 
\end{defn}

The set of 
    smooth $\Gam$-semi-equivariant 
    $(G,\tht)$-valued transition cocycles 
    over $\cU$ 
    will be denoted 
    $\TCocycles_\Gam(\cU,X,(G,\tht))$. 
The set of equivalence classes of
    smooth $\Gam$-semi-equivariant 
    $(G,\tht)$-valued transition cocycles 
    over $\cU$ 
    will be denoted 
    by $\TCocycles_\Gam^\iso(\cU,X,(G,\tht))$.

The first step toward a correspondence 
    between principal bundle and cocycles, 
    is to show how 
    a semi-equivariant transition cocycle 
    can be constructed from 
    a semi-equivariant principal bundle.
    Implicit in the proof of this result 
    is the derivation 
    of the semi-equivariant cocycle property.

\begin{prop}
\label{prop:TC-from-PB}
Let $P \in \PBundles_\Gam(X,(G,\tht))$ and 
        $s := \set{s_a: U_a \rightarrow P|_{U_a}}$
    be a choice of smooth local sections 
    over the cover $\cU$.
    The collection of maps  
    \begin{equation*}
          \phi^{s}
            :=
          \set{
              \phi_{ba}(\gam,\cdot) 
                : 
              U_a \cap \gam^\inv U_b \rightarrow G 
                \mid
              U_a \cap \gam^\inv U_b \neq \emptyset
              } 
    \end{equation*}
    defined by
    \begin{equation}
    \label{eq:DefiningACocycle}
        \gam s_a(x) = s_b(\gam x)\phi_{ba}(\gam,x).
    \end{equation}
    is a smooth $\Gam$-semi-equivariant 
    $(G,\tht)$-valued transition cocycle.
    \begin{proof}
    The given condition 
    implies the following three identities
        \begin{align*}
            \gam'\gam s_a(x) 
                &= 
            s_c(\gam'\gam x)\phi^s_{ca}(\gam'\gam,x)
            &\\
            \gam's_b(\gam x) 
                &= 
            s_c(\gam'\gam x)\phi^s_{cb}(\gam',\gam x)
            &\\
            \gam s_a(x) 
                &= 
            s_b(\gam x)\phi^s_{ba}(\gam,x),
        \end{align*}
    which, together, imply
        \begin{align*}
          s_c(\gam'\gam x)\phi^s_{ca}(\gam'\gam,x) 
            &= \gam'\gam s_a(x) \\
            &= \gam'(s_b(\gam x)\phi^s_{ba}(\gam,x)) \\
            &= (\gam' s_b(\gam x))(\gam'\phi^s_{ba}(\gam,x)) \\
            &= s_c(\gam'\gam x)\phi^s_{cb}(\gam',\gam x)(\gam'\phi^s_{ba}(\gam,x)).
        \end{align*}
    Thus, 
    $\phi^s$ satisfies the cocycle property
    \begin{equation*}
    \phi^s_{ca}(\gam'\gam,x) 
    = 
    \phi^s_{cb}(\gam',\gam x)
    (\gam'\phi^s_{ba}(\gam,x)).
    \end{equation*}
    \end{proof}
\end{prop}
Note that
        \eqref{eq:DefiningACocycle} 
    is the defining relation for a non-equivariant transition cocycle when restricted to $\gam = 1$.
    If $\tht = \id$, then 
        \eqref{eq:DefiningACocycle} 
    is the defining relation 
    for an equivariant transition cocycle.
    
The map
    from 
        semi-equivariant principal bundles to
        semi-equivariant transition cocycles,
    defined by Proposition \ref{prop:TC-from-PB},
    depends on a choice of local sections. 
    However, if one passes to
    isomorphism classes of principal bundles 
    and 
    equivalence classes of transition cocycles
    this dependence disappears.  
    The next proposition shows 
    that cocycles associated to 
    isomorphic principal bundles 
    by Proposition \ref{prop:TC-from-PB}
    are always equivalent, 
    regardless of which sections are chosen.
\begin{prop}
Let
    $P_i \in \PBundles_\Gam(X,(G,\tht))$, and 
    $\phi^i \in \TCocycles_\Gam(\cU,X,(G,\tht))$ 
    be the cocycles associated to local sections
        $s^i := \set{s^i_a: U_a \rightarrow P_i|_{U_a}}$
    as in Proposition \ref{prop:TC-from-PB}.
If $\vphi: P_1 \rightarrow P_2$ is an isomorphism,
  then the collection of maps 
    \begin{equation}
        \mu := \set{\mu_a : U_a \rightarrow G}
    \end{equation}
  defined by 
    \begin{equation}\label{ref:defprop0}
        \vphi(s^1_a(x)) := s^2_a(x)\mu_a(x)
    \end{equation}
    is an equivalence between $\phi^1$ and $\phi^2$.
\begin{proof}
The properties of semi-equivariant principal bundle isomorphisms
    and the defining property \eqref{ref:defprop0} imply that
    \begin{align*}
      \vphi(\gam s^1_a(x)) &= \gam \vphi(s^1_a(x)) \\
      \vphi( s^1_b(\gam x) \phi^1_{ba}(\gam,x) ) &= \gam( s^2_a(x)\mu_a(x) )  \\
      \vphi( s^1_b(\gam x)) \phi^1_{ba}(\gam,x) &= (\gam s^2_a(x)) (\gam \mu_a(x) )  \\
      s^2_b(\gam x)\mu_b(\gam x) \phi^1_{ba}(\gam,x) &= s^2_b(\gam x)\phi^2_{ba}(\gam,x) (\gam \mu_a(x) ).
    \end{align*}
Thus,
    \begin{equation*}
        \mu_b(\gam x) \phi^1_{ba}(\gam,x) = \phi^2_{ba}(\gam,x) (\gam \mu_a(x) ),
    \end{equation*}
    and $\mu$ is an equivalence between  $\phi^1$ and $\phi^2$ for any choice of sections $s^i$.
\end{proof}
\end{prop}

\begin{cor}
The map of Proposition \ref{prop:TC-from-PB} induces a well-defined map
    \begin{align*}
          \PBundles^{\iso}_\Gam(X, (G,\tht)) &\rightarrow \TCocycles^{\iso}_\Gam(\cU, X, (G,\tht))  \\
            [P] &\mapsto [\phi^s],
    \end{align*}
    where $s$ is 
    any collection of smooth local sections 
    of $P$.
\end{cor}

The correspondence between semi-equivariant cocycles and principal bundles 
    has now been shown in one direction.
    Next, an inverse map reconstructing a semi-equivariant principal bundle
    from a semi-equivariant transition cocycle is defined.
\begin{prop}\label{prop:PB-from-TC}
  Let $\phi \in \TCocycles_\Gam(\cU, X, (G,\tht))$.
    The bundle $P^\phi$ defined by
    \begin{equation*}
       \pi: ( \dUnion_{a \in A} U_a \times G / \sim ) \rightarrow X, 
    \end{equation*}
  where
  \begin{enumerate}
    \item $(a,x,g) \sim (b,x,\phi_{ba}(1,x)g)$ defines the equivalence relation $\sim$
    \item $\pi[a,x,g] := x$ is the projection map 
    \item $[a,x,g]g' := [a,x,gg']$ defines the right-action of $G$
    \item $\gam[a,x,g] := [b,\gam x, \phi_{ba}(\gam,x) (\gam g)]$ defines the left action of $\Gam$,
  \end{enumerate}
  is a 
  smooth 
  $\Gam$-semi-equivariant 
  principal $(G,\tht)$-bundle.
\begin{proof}
The elements $\set{\phi_{ba}(1, \cdot)}$ satisfy
    \begin{equation*}
        \phi_{ca}(1, x) = \phi_{cb}(1, x)\phi_{ba}(1, x)
    \end{equation*}
  and so form a $G$-valued cocycle in the usual sense.
  Therefore, the usual proof that $P^\phi$ is a principal $G$-bundle applies.
The $\Gam$-action is well-defined on equivalence classes 
  as
  \begin{align*}
    \gam[b,x,\phi_{ba}(1,x)g] 
      &= [c,\gam x, \phi_{cb}(\gam,x) \gam(\phi_{ba}(1,x)g)] \\
      &= [c,\gam x, \phi_{cb}(\gam,x) (\gam\phi_{ba}(1,x)) (\gam g)] \\
      &= [c,\gam x, \phi_{ca}(\gam,x) (\gam g)] \\
      &= \eta_\gam[a,x,g].
  \end{align*}
The semi-equivariance property
    $\gam(pg) = (\gam p)(\gam g)$ 
  is satisfied as
  \begin{align*}
    \gam([a,x,g]g')
      &=
    \gam([a,x,gg']) \\
      &=
    [b,\gam x, \phi_{ba}(\gam,x) (\gam g g')] \\
      &=
    [b,\gam x, \phi_{ba}(\gam,x) (\gam g) (\gam g')] \\
      &=
    (\gam[a,x,g]) (\gam g')
  \end{align*}
Thus, $P^\phi$ is a $\Gam$-semi-equivariant principal $(G,\tht)$-bundle.
\end{proof}
\end{prop}

This reconstruction map is also well-defined at the level of
isomorphism and equivalence classes.

\begin{prop}
Let $\phi^i \in \TCocycles_\Gam(\cU,X,(G,\tht))$ and $P_i \in \PBundles_\Gam(X,(G,\tht))$
    be the associated principal bundles, 
    constructed using Proposition \ref{prop:PB-from-TC}.
    If $\mu := \set{\mu_a : U_a \rightarrow G}$
    is an equivalence between $\phi^1$ and $\phi^2$
    then
    \begin{align*}
      \vphi: P_1 &\rightarrow P_2 \\
      [a,x,g] &\mapsto [a,x,\mu_a(x)g].
    \end{align*}
    is an isomorphism.
\begin{proof}
    That $\vphi$ is a well-defined isomorphism of principal $G$-bundles
    follows immediately from the proof in the non-equivarant case. 
    Compatibility with the $\Gam$-action is satisfied as
        \begin{align*}
          \gam \vphi([a,x,g])
            &=
          \gam[a,x,\mu_a(x)g] \\
            &=
          [b,\gam x,\phi'_{ba}(\gam,x) \gam(\mu_a(x)g)] \\
            &=
          [b,\gam x,\phi'_{ba}(\gam,x) (\gam\mu_a(x))(\gam g)] \\
            &=
          [b,\gam x, \mu_b(\gam x)\phi_{ba}(\gam,x) (\gam g)] \\
            &=
          \vphi([b,\gam x, \phi_{ba}(\gam,x) (\gam g)]) \\
            &=
          \vphi( \gam [a,x,g] ).
        \end{align*}
    Thus, $\vphi$ is an isomorphism of $\Gam$-semi-equivariant principal $(G,\tht)$-bundles.
\end{proof}
\end{prop}

\begin{cor}
The map of Proposition \ref{prop:PB-from-TC} induces a well-defined map
    \begin{align}
           \TCocycles^{\iso}_\Gam(\cU, X, (G,\tht)) &\rightarrow \PBundles^{\iso}_\Gam(X, (G,\tht))  \\
            [\phi] &\mapsto [P^\phi].
    \end{align}
\end{cor}

Finally, one shows that the two maps defined above are inverse to one another.

\begin{prop}
The maps
    \begin{align*}
           \TCocycles^{\iso}_\Gam(\cU, X, (G,\tht)) &\rightarrow \PBundles^{\iso}_\Gam(X, (G,\tht))  \\
            [\phi] &\mapsto [P^\phi]
    \end{align*}
    and
    \begin{align*}
          \PBundles^{\iso}_\Gam(X, (G,\tht)) &\rightarrow \TCocycles^{\iso}_\Gam(\cU, X, (G,\tht))  \\
            [P] &\mapsto [\phi^s]
    \end{align*}
    are inverse to one another.
\begin{proof}
Let $P \in \PBundles_\Gam(X,(G,\tht))$,
    $\phi := \phi^s$ and $P' := P^\phi$
  for some collection of local sections 
    $s := \set{s_a: U_a \rightarrow P|_{U_a}}$. 
The sections $\set{s_a}$ define a trivialization 
    $\set{t_a}$ of $P$ 
  by
  \begin{align*}
    t_a: P|_{U_a} &\rightarrow U_a \times G \\
    s_a(x) &\mapsto (a,x,1)
  \end{align*}
  and a collection of maps
    $\set{T_a: P|_{U_a} \rightarrow G}$ 
  by
    $t_a(p) =: (a,x,T_a(p))$
  where $x = \pi_P(p)$.
  Note that 
    $T_a(pg) = T_a(p)g$.
Define 
  \begin{align*}
    \vphi: P &\rightarrow P'  \\
        p &\mapsto [t_a(p)].
  \end{align*}
  That $\vphi$ is a well-defined isomorphism of
  principal $G$-bundles follows from the proof
  in the non-equivariant case.
To check that $\vphi$ is compatible with the $\Gam$-actions
  first note that
  \begin{align*}
    t_b \circ \eta_\gam \circ t^\inv_a(a,x,g)
      &=
    t_b(\gam(s_a(x)g)) \\
      &=
    t_b((\gam s_a(x)) (\gam g)) \\
      &=
    t_b(s_b(\gam x)\phi_{ba}(\gam,x) (\gam g)) \\
      &=
    (b, \gam x, \phi_{ba}(\gam,x) (\gam g))
  \end{align*}
  where $\eta$ is the $\Gam$-action on $P$.
  Thus,
  \begin{align*}
    \gam \vphi(p)
      &=
    \gam[t_a(p)] \\
      &=
    \gam[a,x,T_a(p)] \\
      &=
    [b,\gam x,\phi_{ba}(\gam,x)T_a(p)] \\
      &=
    [t_b \circ \eta_\gam \circ t^\inv_a(a,x,T_a(p))] \\
      &=
    [t_b(\gam p)] \\
      &=
    \vphi(\gam p).
  \end{align*}
Therefore, $\vphi$ is an isomorphism of $\Gam$-semi-equivariant principal $(G,\tht)$-bundles and 
      $P \mapsto \phi^s \mapsto P^{\phi^s}$  
    is the identity map at the level of isomorphism classes.
\end{proof}
\end{prop}

The main theorem of this section has now been proved.

\begin{thm}
\label{thm:PB-TC-bij}
There is a bijective correspondence
    \begin{equation*}
    \PBundles^{\iso}_\Gam(X,(G,\tht)) 
    \leftrightarrow 
    \TCocycles^{\iso}_\Gam(\cU,X, (G,\tht))
    \end{equation*}
    between semi-equivariant cocycles and principal bundles.
\end{thm}

\pagebreak
\begin{center}
    \tikz[scale=1]{\tikzstyle{every node}=[font=\tiny]
        \draw [->>] (0,4) -- (0,6);
        \draw [->] (-0.5,4.5) -- (0.5,5.5);

        \draw [->] (1.75,5) arc [start angle=340, end angle=20, radius=0.25cm]; 
        \node at (1.5,5) [above=14,align=center] {clockwise\\rotation};

        \draw [<-] (3,4) -- (3,6);
        \draw [->>] (2.5,4.5) -- (3.5,5.5);

        \draw [->] (4.5,4.75) arc [start angle=-90, end angle=90, radius=0.25cm]; 
        \draw (4.25,4.75) -- (4.75,5.25);
        \node at (4.5,4.75) [above=27,align=center] {reflection};

        \draw [->] (6,4) -- (6,6);
        \draw [->>] (5.5,4.5) -- (6.5,5.5);

        \draw [->,color=blue] (1.5,4.53)
                     .. controls +(down:5mm) and +(up:5mm) .. 
                   (4.5,2.23);

        \draw [->>] (0,0) -- (0,2);
        \draw [->] (-0.5,0.5) -- (0.5,1.5);

        \draw [->] (1.5,0.75) arc [start angle=-90, end angle=90, radius=0.25cm]; 
        \draw (1.25,0.75) -- (1.75,1.25);
        \node at (1.5,0.75) [above=27,align=center] {reflection};

        \draw [<<-] (3,0) -- (3,2);
        \draw [->] (2.5,0.5) -- (3.5,1.5);

        \draw [<-] (4.75,1) arc [start angle=340, end angle=20, radius=0.25cm];
        \node at (4.5,1) [above=14,align=center] {anti-clockwise\\rotation};

        \draw [->] (6,0) -- (6,2);
        \draw [->>] (5.5,0.5) -- (6.5,1.5);

        \node at (6,3) [align=center,font=\small] {$=$};
        \node at (0,3) [align=center,font=\small] {$=$};
        }


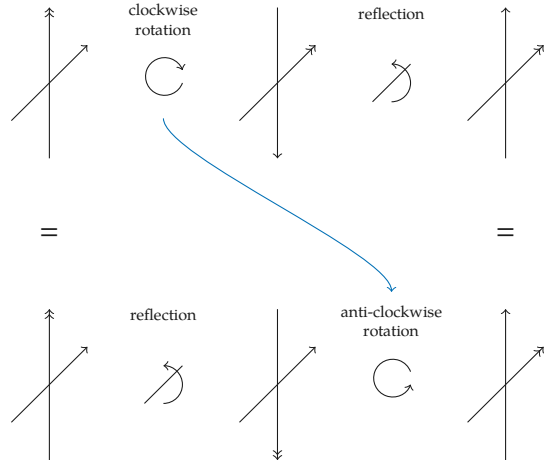
\captionof{figure}{
            This figure corresponds to $\C$ 
            equipped with 
            conjugation as a $\Z_2$-action and 
            $\Unitary(1)$ acting by rotations,
            as in Example \ref{exam:UKSDP}.
            The blue line represents 
            the conjugation automorphism 
            on $\Unitary(1)$.
            This conjugation 
            is required 
            in order to obtain 
            the same final result
            when the two actions are applied 
            in reversed order.
            } \label{fig:SEexam}
\end{center}

\begin{center}

    \tikz[scale=1.2]{\tikzstyle{every node}=[font=\tiny]
       \draw (0,1) ellipse [x radius=1cm, y radius=0.4cm] node[below right=11pt] {$U_a$}; 
       \node[circle,fill=black,inner sep=0pt,minimum size=1.25pt,label=right:{$x$}] at (0,1) {};
       \node(Sa) at (0,2) {$s_a$};

       \draw (3,1) ellipse [x radius=1cm, y radius=0.4cm] node[below right=11pt] {$U_b$}; 
       \node[circle,fill=black,inner sep=0pt,minimum size=1.25pt,label=right:{$\gam x$}] at (3,1) {};
       \node(Sb)  at (3,2) {$s_b$};
       \node(gSa)  at (3,4) {$\gam s_a$}
            edge [<-] node[auto,swap] {$\gam$} (Sa)
            edge [<-] node[auto,swap] {$\phi_{ba}(\gam,x)$} (Sb);
 
       \draw (6,1) ellipse [x radius=1cm, y radius=0.4cm] node[below right=11pt] {$U_c$}; 
       \node[circle,fill=black,inner sep=0pt,minimum size=1.25pt,label=right:{$\gam'\gam x$}] at (6,1) {};
       \node(Sc) at (6,2) {$s_c$};
       \node(gSb) at (6,4) {$\gam' s_b$}
            edge [<-] node[auto, swap] {$\gam'$}(Sb)
            edge [<-] node[auto] {$\phi_{cb}(\gam',\gam x)$}(Sc);

       \node(ggSa) at (6,6) {$\gam'\gam s_c$}
            edge [<-] node[auto,swap] {$\gam'$}(gSa)
            edge [<-, dashed] node[auto] {$\mathbin{\color{blue}\gam'}\phi_{ba}(\gam, x)$}(gSb);
       
        \draw [->] 
            (Sc.east) 
                .. controls +(right:25mm) and +(right:25mm) .. 
            (ggSa.east)
            node[midway,anchor=west] {$\phi_{ca}(\gam'\gam, x)$}; {};
    
    \node(Comm) [align=center] at (4.5,4) {Semi-\\equivariance}; 
    \node(Comm) [align=center] at (7.5,4) {Cocycle\\property}; 
    }
    
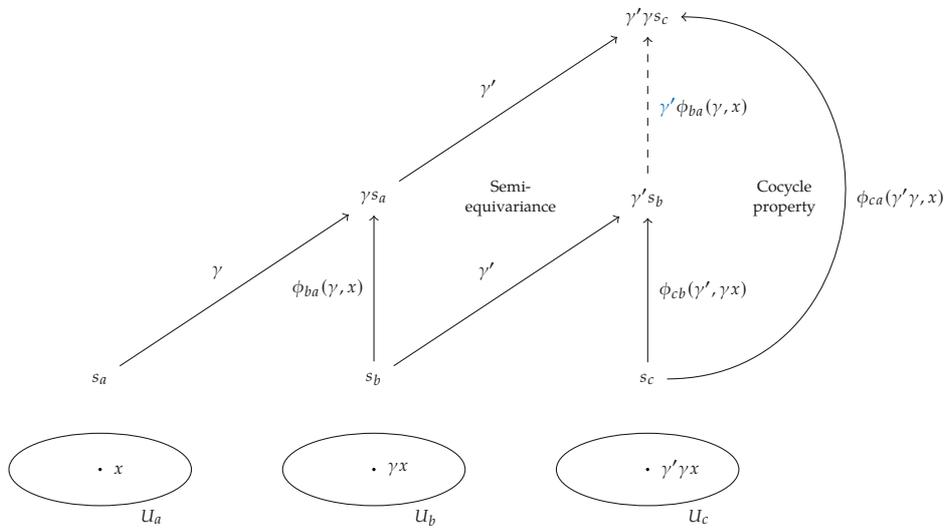
\captionof{figure}{
        This diagram represents the derivation 
        of the semi-equivariant cocycle property,
        as in Proposition \ref{prop:TC-from-PB}.
        Each node of the diagram represents 
        a local section of a principal bundle.
        The diagonal arrows represent applications 
        of the $\Gam$-action,
        while the vertical arrows represent 
        the action of a cocycle $\phi$
        via the right action
        of the structure group.
        With the exception of the dashed line,
        all of the arrows follow from the definitions.
        The dashed line follows 
        by the semi-equivariance property 
        of the principal bundle,
        the blue $\gam'$ is acting 
        on the element $\phi_{ba}(\gam,x)$ 
        of the structure group.
    } \label{fig:SECprop}
\end{center}

\begin{figure}[H]
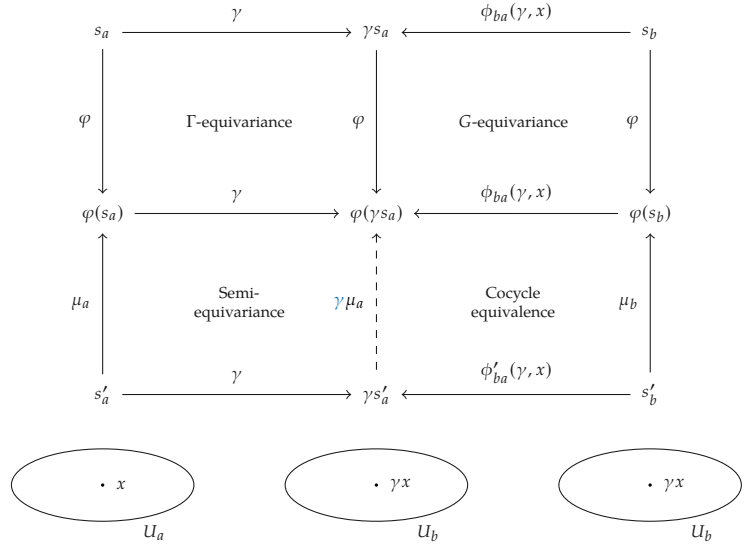

    \centering
\tikz[scale=1.2]{\tikzstyle{every node}=[font=\tiny]
   \draw (0,1) ellipse [x radius=1cm, y radius=0.4cm] node[below right=11pt] {$U_a$}; 
   \node[circle,fill=black,inner sep=0pt,minimum size=1.25pt,label=right:{$x$}] at (0,1) {};
   \node(S2a)    at (0,2) {$s'_a$};
   \node(Sa)     at (0,6) {$s_a$};

   \draw (3,1) ellipse [x radius=1cm, y radius=0.4cm] node[below right=11pt] {$U_b$}; 
   \node[circle,fill=black,inner sep=0pt,minimum size=1.25pt,label=right:{$\gam x$}] at (3,1) {};

   \draw (6,1) ellipse [x radius=1cm, y radius=0.4cm] node[below right=11pt] {$U_b$}; 
   \node[circle,fill=black,inner sep=0pt,minimum size=1.25pt,label=right:{$\gam x$}] at (6,1) {};
   \node(S2b)    at (6,2) {$s'_b$};
   \node(Sb)     at (6,6) {$s_b$};

   \node(gSa)     at (3,6) {$\gam s_a$}
        edge [<-] node[auto,swap] {$\gam$} (Sa)
        edge [<-] node[auto] {$\phi_{ba}(\gam,x)$} (Sb);
   \node(mSa)            at (0,4) {$\vphi(s_a)$}
        edge [<-] node[auto] {$\vphi$} (Sa)
        edge [<-] node[auto,swap] {$\mu_a$} (S2a);
   \node(gS2a)                 at (3,2) {$\gam s'_a$}
        edge [<-] node[auto,swap] {$\gam$} (S2a)
        edge [<-] node[auto] {$\phi'_{ba}(\gam,x)$} (S2b);
   \node(mSb)                  at (6,4) {$\vphi(s_b)$}
        edge [<-] node[auto] {$\vphi$} (Sb)
        edge [<-] node[auto,swap] {$\mu_b$} (S2b);

   \node(mgSa)     at (3,4) {$\vphi(\gam s_a)$}
        edge [<-,dashed] node[auto,swap] {$\mathbin{\color{blue}\gam}\mu_a$} (gS2a)
        edge [<-] node[auto] {$\vphi$} (gSa)
        edge [<-] node[auto,swap] {$\gam$} (mSa)
        edge [<-] node[auto] {$\phi_{ba}(\gam,x)$} (mSb);

    \node(GamEqui) at (1.5,5) {$\Gam$-equivariance}; 
    \node(GEqui) at (4.5,5) {$G$-equivariance}; 
    \node(Comm) [align=center] at (1.5,3) {Semi-\\equivariance}; 
    \node(Equiv) [align=center] at (4.5,3) {Cocycle\\equivalence}; 
}
    \caption{
        This diagram depicts the derivation 
        of the equivalence property 
        for semi-equivariant cocycles,
        see Definition \ref{defn:SETCEquiv}.
        Here $\vphi$ is 
        a semi-equivariant principal bundle isomorphism.
        Each node of the diagram 
        represents a local section 
        of a principal bundle.
        The arrows running downward 
        are applications of 
        a principal bundle isomorphism $\vphi$.
        The arrows running left to right 
        are applications of the $\Gam$-action.
        The arrows running right to left 
        are right actions by the cocycle $\phi$.
        Those running upward are
        right actions of the cocycle equivalence $\mu$.
        With the exception of the dashed arrow, 
        all of the arrows follow from definitions.
        The commutation of the top two squares follows 
        from the properties 
        of principal bundle isomorphisms.
        The dashed arrow is follows 
        from the semi-equivariance property 
        of the principal bundle.
        This twists the equivalence $\mu_a$ 
        by the action of $\Gam$ 
        on the structure group, which is marked in blue.
        The lower right square is 
        the semi-equivariant cocycle equivalence condition.
    }\label{fig:SEEprop}
\end{figure}
\pagebreak


\section{Semi-equivariant Cohomology}
\label{sec:semiCo}

In order to study liftings 
    of semi-equivariant principal bundles, 
    a cohomology theory is needed.
    The existing notions 
    of equivariant cohomology 
    are inappropriate for this task, 
    and a new cohomology theory must be constructed.
In this section, 
    a $\Gam$-semi-equivariant 
    \v{C}ech cohomology theory 
    is developed
    with an abelian $\Gam$-group $(G,\tht)$ 
    as its coefficient group.
    The theory makes use of a simplicial space 
    which 
    encodes the group structure of $\Gam$, and 
    the action of $\Gam$ 
    on the manifold $X$.
    In addition to these actions, 
    the effect of the action $\tht$ 
    must be incorporated.
    This is achieved 
    by twisting the coboundary map using $\tht$.
    There are a few details to check, 
    but everything works as one would wish.
    This semi-equivariant cohomology theory generalises 
    an equivariant cohomology theory 
    outlined by Brylinski
        \cite[\S A]{Brylinski.2000.Gerbes-on-complex-reductive-Lie-groups}
    Another helpful reference is 
    \cite[\S3.3]{Gomi.2005.Equivariant-smooth-Deligne-cohomology}.
    One feature of the presentation here 
    is that it avoids the use of hypercohomology.
    The second dimension of the bicomplex appearing in 
        \cite[\S A]{Brylinski.2000.Gerbes-on-complex-reductive-Lie-groups}
    is an artifact of the choice 
    to separate the cocycle into two parts, 
    one encoding the transition functions 
    for the total space and 
    one encoding the action.
    Although this is ultimately a notational matter, 
    the reduced book-keeping is helpful 
    when checking higher cocycle conditions.

The construction of 
    semi-equivariant \v{C}ech cohomology 
    begins with the definition
    of a simplicial space.
    The coboundary map on 
    the underlying chain complex 
    of the cohomology theory
    will be constructed using
    the face maps of this space.
\begin{defn}
  \label{def:}
  Let $X$ be a manifold equipped 
  with a smooth action of $\Gam$.
The \emph{simplicial space} associated to $X$
  is defined by
    \[ X^\bul := \set{\Gam^p \times X}_{p \geq 0}. \]
The simplicial space carries 
  \emph{face} and \emph{degeneracy} maps
  \begin{align*}
       d^{p}_i : X^p &\rightarrow X^{p-1} &
       e^{p}_i : X^p &\rightarrow X^{p+1}
  \end{align*}  
  defined by
  \begin{align}
      d^{p}_i(\gam_1,\ldots,\gam_p,x)
        &:=
      \begin{cases}
        (\gam_2, \ldots, \gam_p,x) 
          &\txt{ for } i=0 \\
        (\gam_1, \ldots,\gam_i\gam_{i+1},\ldots,\gam_p,x) 
          &\txt{ for } 1 \leq i \leq p-1 \\
        (\gam_1, \ldots,\gam_{p-1},\gam_p x) 
          &\txt{ for } i=p \\
      \end{cases} \label{al:SimpSpaceFaceMaps}
    \\
      e^{p}_i(\gam_1,\ldots,\gam_{p},x)
        &:=
        (\gam_1,\ldots,\gam_i,1,\gam_{i+1},\ldots,\gam_p,x) 
        \phantom{..} \txt{ for } 0 \leq i \leq p+1 \nonumber
  \end{align}
\end{defn}
Notice that 
    in \eqref{al:SimpSpaceFaceMaps} 
    the face map $d^p_0$ 
    discards the element $\gam_1$,
    this element will be used 
    to define the simplicial twisting maps, 
    in Definition \ref{def:SimpTwist}.
\begin{prop}{}
  \label{prop:}
The face and degeneracy maps satisfy 
  the simplicial identities
\begin{equation}
  \begin{aligned}
        d_i \circ d_j 
            &=
        d_{j-1} \circ d_i &\txt{ for } i < j     
        \\
        e_i \circ e_j 
            &= 
        e_{j+1} \circ e_i &\txt{ for } i \leq j 
    \end{aligned}
    \hspace{1.2cm}
    \begin{aligned}
        d_i \circ e_j &=
        \begin{cases}
          e_{j-1} \circ d_i &\txt{ for } i < j          \\
          \id               &\txt{ for } i = j, j+1     \\
          e_j \circ d_{i-1} &\txt{ for } i > j+1   
        \end{cases}    
   \end{aligned}
\end{equation}
\end{prop}
Corresponding to the face maps 
    $d^p_i$, 
  twisting maps
    $\tht_i: X^p \times G \rightarrow G$
  can be defined.
  These maps encode 
    the action $\tht$ of $\Gam$ on $G$
    and will be used to twist 
    the coboundary map.
    They are the basic ingredient needed 
    for generalisation 
    to the semi-equivariant setting.
    Note that 
    it is only the twisting map $\tht_0$ 
    that has any effect. 
    The rest of the twisting maps 
    are included for notational convenience 
    when dealing with simplical identities.
\begin{defn}{}
  \label{def:SimpTwist}
The \emph{simplicial twisting maps} 
    $\tht_i: X^p \times G \rightarrow G$
  are given by
\begin{align*}
    \tht^{(\gam_1, \ldots, \gam_{p},x)}_i
      &:=
    \begin{cases}
       \tht_{\gam_1}
        &\txt{ for } i=0 \\
      \id
        &\txt{ for } 1 \leq i \leq p-1 \\
      \id
        &\txt{ for } i=p \\
    \end{cases}
\end{align*}
\end{defn}

The twisting maps also satisfy simplicial identities
    which help to ensure that the coboundary 
    map in semi-equivariant cohomology squares to zero.

\begin{prop}{}
  \label{prop:}
  The simplicial twisting maps satisfy the identities
  \begin{align*}
    \tht^{x^{p+1}}_j \circ \tht^{d_j(x^{p+1})}_i
      &=
    \tht^{x^{p+1}}_i \circ \tht^{d_i(x^{p+1})}_{j-1}
    \hspace{0.5cm}
    \txt{ for } i < j    \\
    \tht^{e_j(x^p)}_i &= 
    \begin{cases}
    \tht_i^{x^p} &\txt{ for } i < j    \\
    \id & \txt{ for } i = j,j+1    \\
    \tht^{x^p}_{i-1} &\txt{ for } i > j+1,
    \end{cases}
  \end{align*}
    where $x^p \in X^p$.
\begin{proof}
  The identities are trivial 
  for most combinations of $i$ and $j$.
  The remaining cases can be checked individually.
  In particular, the first identity reduces to
    \begin{equation*}
        \begin{aligned}
            \id \circ \tht_{\gam_1\gam_2} 
                &= 
            \tht_{\gam_1} \circ \tht_{\gam_2}
            \\
            \id \circ \tht_{\gam_1} 
                &= 
            \tht_{\gam_1} \circ \id
            \\
            \id &= \id 
        \end{aligned}
        \hspace{1cm}
        \begin{aligned}
            &\txt{ for } i=0,j=1 \\
            &\txt{ for } i=0,j \geq 2 \\
            &\txt{ otherwise.}
        \end{aligned}
    \end{equation*}
\end{proof}
\end{prop}
To construct a \v{C}ech-type theory,
    a simplicial cover 
        $\cU^\bul$ of $X^\bul$ 
    is needed.
    Such a cover can be constructed 
    from an appropriate cover 
        $\cU := \set{U_a \mid a \in A}$ of $X$.
    First, 
    the indexing set of the simplicial cover 
    is defined.
    This indexing set has a simplicial structure 
    defined by face and degeneracy maps,
    which will again be denoted by $d^p_i$ and $e^p_i$.
\begin{defn}{}
  \label{def:}
Define the indexing set for $\cU^\bul$ by
    \[ A^\bul := \set{A^p}_{p \geq 0} \]
  where
    $A^p := \set{ (a_0, \ldots, a_p) \mid a_i \in A }$.
  Elements of $A^p$ will be denoted by $a^p$.
  This set carries face and degeneracy maps
  \begin{align*}
      d^{p}_i : A^p &\rightarrow A^{p-1} &
      e^{p}_i : A^p &\rightarrow A^{p+1}
  \end{align*}
  defined by
  \begin{align*}
    d^{p}_i(a_0,\ldots,a_p) 
      &:= 
    (a_0, \ldots, \hat{a}_i, \ldots, a_p) \\
    e^{p}_i(a_0,\ldots,a_p) 
      &:= 
    (a_0, \ldots, a_i, a_i, a_{i+1}, \ldots, a_p),
  \end{align*}
    where $\hat{a}_i$ denotes 
        the removal of the element $a_i$.
\end{defn}
\begin{prop}{}
  \label{prop:}
The face and degeneracy maps of the indexing set $A^\bul$ satisfy 
    \begin{equation*}
        \begin{aligned}
            d_i \circ d_j &= d_{j-1} \circ d_i &\txt{ for } i < j     \\  
            e_i \circ e_j &= e_{j+1} \circ e_i &\txt{ for } i \leq j
        \end{aligned}
        \hspace{1.2cm}
        \begin{aligned}
            d_i \circ e_j =
            \begin{cases}
              e_{j-1} \circ d_i &\txt{ for } i < j          \\
              \id               &\txt{ for } i = j, j+1     \\
              e_j \circ d_{i-1} &\txt{ for } i > j+1.   
            \end{cases}    
        \end{aligned}
    \end{equation*}
\end{prop}
Before defining the simplicial cover itself, observe that 
    the elements of the simplicial space define sequences of points in $X$.
\begin{defn}
    Let 
        $x^p = (\gam_1,\ldots,\gam_p, x) \in X^p$.
    The associated sequence 
        $\set{\mb{x}^p_i}$ 
    is defined by
    \begin{equation*}
        \mb{x}^p_i := \gam_{p-i} \cdots \gam_p x \in X.
    \end{equation*}
\end{defn}

Simplicial covers generalise the nerves of covers.
    The definition will be made using the definitions of
    the sequences $\mb{x}^p_i$ and indexing set $A^\bul$.
\begin{defn}
The \emph{simplicial cover} 
    \begin{equation*}
        \cU^\bul := \set{\cU^p}_{p \geq 0}
    \end{equation*}
    associated to $\cU$ is a sequence of covers
    $\cU^p$ of $X^p$ each indexed by $A^p$.
    A set 
    \begin{equation*}
        U_{(a_0,\ldots,a_p)} \in \cU^p
    \end{equation*}
    consists of all points in $X^p$ 
    such that $\mb{x}^p_i \in U_{a_i}$ 
    for $0 \leq i \leq p$.
\end{defn}
For example, 
        $
            (\gam_1,\gam_2,\gam_3,x) 
                \in 
            U_{(a_{0},a_{1},a_{2},a_{3})}
        $
    can be visualised as a path
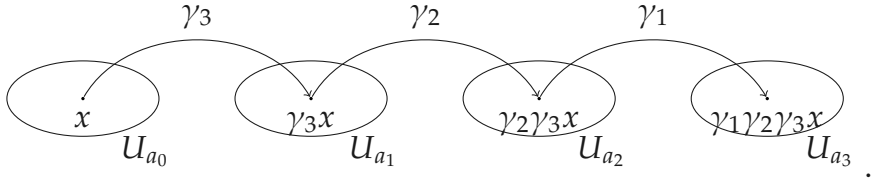
\begin{center}
    \tikz[scale=1]{
       \draw (0,1) ellipse [x radius=1cm, y radius=0.5cm] node[below right=10pt] {$U_{a_0}$}; 
       \node(x)[circle,fill=black,inner sep=0pt,minimum size=1.25pt,label=below:{$x$}] at (0,1) {};

       \draw (3,1) ellipse [x radius=1cm, y radius=0.5cm] node[below right=10pt] {$U_{a_1}$}; 
       \node(gx)[circle,fill=black,inner sep=0pt,minimum size=1.25pt,label=below:{$\gam_3 x$}] at (3,1) {}
            edge[<-,bend right=60] node [midway,above] {$\gam_3$} (x);
 
       \draw (6,1) ellipse [x radius=1cm, y radius=0.5cm] node[below right=10pt] {$U_{a_2}$}; 
       \node(ggx)[circle,fill=black,inner sep=0pt,minimum size=1.25pt,label=below:{$\gam_2\gam_3 x$}] at (6,1) {}
            edge[<-,bend right=60] node [midway,above] {$\gam_2$} (gx);

       \draw (9,1) ellipse [x radius=1cm, y radius=0.5cm] node[below right=10pt] {$U_{a_3}$}; 
       \node(gggx)[circle,fill=black,inner sep=0pt,minimum size=1.25pt,label=below:{$\gam_1\gam_2\gam_3 x$}] at (9,1) {}
            edge[<-,bend right=60] node [midway,above] {$\gam_1$} (ggx);
    }.
\end{center}
Note that a refinement of $\cU$ induces a refinement of $\cU^\bul$.
    Also, the face maps of the simplicial cover are compatible
    with those of the simplicial space.
    This is necessary to ensure that the coboundary map is well-defined.
\begin{prop}
    The pullback maps of the simplicial space are compatible with
    those on the indexing set of the cover in the sense that
        $d_i(U_{a^p}) \subseteq U_{d_i(a^p)}$.
\end{prop}


Semi-equivariant \v{C}ech cohomology 
    is based on a single cochain complex.
    A $p$-cochain for this cohomology theory 
    consists of a smooth function 
    on each set in the $p$th level 
    of the simplicial cover.
\begin{defn}\label{def:cochains}
The group of \emph{$p$-cochains} 
    is defined by
    \begin{equation*}
        K_\Gam^p(\cU,X,(G,\tht)) 
          := 
        \prod_{a^p \in A^p} C^\infty(U_{a^p}, G),
    \end{equation*}
  with the group operation
    $(\phi'\phi)_{a^p} := \phi'_{a^p}\phi_{a^p}$.
\end{defn}
These cochains can be pulled back 
    by the face maps.
    In the semi-equivariant setting, 
    the pullback maps are composed with the twisting maps.
    This modifies the pullback by $d_0$.
\begin{defn}\label{def:TPB}
  The \emph{twisted pullback maps}
    \begin{equation*}
        \d^p_i: K_\Gam^p(\cU,X,(G,\tht)) \rightarrow K_\Gam^{p+1}(\cU,X,(G,\tht))
    \end{equation*}
  are defined by
    \begin{equation*}
      (\d^p_i \phi)_{a^{p+1}}(x^{p+1}) 
        := 
      \tht_i^{x^{p+1}} 
        \circ 
      \phi_{d^p_i(a^{p+1})} 
        \circ 
      d^p_i(x^{p+1})
    \end{equation*}
\end{defn}
Note that the property 
    $d_i(U_{a^p}) \subseteq U_{d_i(a^p)}$
  of the cover
  ensures that 
    $\d_i(\phi)$ 
  is a well-defined element of $K_\Gam^{p+1}(\cU,X,(G,\tht))$.
\begin{prop}{}
\label{prop:}
The twisted pullback maps are group homomorphisms.
\begin{proof}
Using the fact that 
    $\tht_\gam$ is an automorphism 
    for all $\gam \in \Gam$,
    \begin{align*}
    (\d_i (\phi'&\phi))_{a^{p+1}}(x^{p+1}) \\
      &= 
    \tht_i^{x^{p+1}} 
      \circ 
    (\phi'\phi)_{d_i(a^{p+1})} 
      \circ 
    d_i(x^{p+1}) \\
      &= 
    \tht_i^{x^{p+1}}
    \Big(
         (\phi'_{d_i(a^{p+1})}\circ d_i(x^{p+1}) )
         (\phi_{d_i(a^{p+1})}\circ d_i(x^{p+1})  )
    \Big) \\
      &= 
        \Big( \tht_i^{x^{p+1}} \circ 
          \phi'_{d_i(a^{p+1})} \circ 
          d_i(x^{p+1}) \Big) 
        \Big(\tht_i^{x^{p+1}} \circ 
          \phi_{d_i(a^{p+1})} \circ 
          d_i(x^{p+1}) \Big) \\
      &=
        \Big( (\d_i \phi')_{a^{p+1}}(x^{p+1}) \Big)
        \Big( (\d_i \phi)_{a^{p+1}}(x^{p+1}) \Big)
    \end{align*}
\end{proof}
\end{prop}
The simplicial identities of the face maps 
    for the simplicial space, the simplicial cover and the twisting maps combine
    to produce a simplicial identity for the twisted pullback maps.
\begin{prop}\label{prop:tpb-si}
For $i<j$ the twisted pullback maps satisfy the identity
    \begin{equation*}
        \d_j \circ \d_i = \d_i \circ \d_{j-1}.
    \end{equation*}
\begin{proof}
Using the corresponding simplicial identities 
  between face maps on the simplicial complex,
  those on the simplicial cover, and
  those between the simplical twisting maps
  one can directly compute
  \begin{align*}
    (\d_j(\d_i \phi))_{a^{p+2}}(x^{p+2})
      &=
    \tht_j^{x^{p+2}} 
      \circ 
    (\d_i \phi)_{d_j(a^{p+2})} 
      \circ 
    d_j(x^{p+2}) \\
      &= 
    \tht_j^{x^{p+2}} 
      \circ 
    \tht_i^{d_j(x^{p+2})} 
        \circ 
    \phi_{d_i \circ d_j(a^{p+2})} 
        \circ 
    d_i
      \circ 
    d_j(x^{p+2}) \\
      &=
    \tht_i^{x^{p+2}} 
      \circ 
    \tht_{j-1}^{d_i(x^{p+2})} 
        \circ 
    \phi_{d_{j-1} \circ d_i(a^{p+2})} 
        \circ 
    d_{j-1}
      \circ 
    d_i(x^{p+2}) \\
      &=
    \tht_i^{x^{p+2}} 
      \circ 
    (\d_{j-1}\phi)_{d_i(a^{p+2})}
      \circ 
    d_i(x^{p+2}) \\
      &=
    (\d_i(\d_{j-1}\phi))_{a^{p+2}}(x^{p+2}).
  \end{align*}
\end{proof}  
\end{prop}
Finally, the coboundary maps are defined.
\begin{defn}\label{def:BM}
The \emph{coboundary maps}
    \begin{equation*}
        \d^p: K_\Gam^p(\cU,X,(G,\tht)) \rightarrow K_\Gam^{p+1}(\cU,X,(G,\tht))
    \end{equation*}
    are defined by
    \begin{equation*}
        \d^p := \sum_{0 \leq i \leq p} (-1)^i \d^p_i.
    \end{equation*}
\end{defn}
Using the simplicial identity 
    for the twisted pullback maps,
    the square of the coboundary map 
    is shown to be zero.
\begin{prop}\label{prop:BMBMeq0}
  The coboundary map satisfies $\d\d = 0$.
\begin{proof}
First note, 
  using Proposition \ref{prop:tpb-si}, 
  that
  \begin{align*}
    \sum_{i < j, j \leq p+2}
      (-1)^{i+j} \d_j\d_i 
      &=
    \sum_{i < j, j \leq p+2}
      (-1)^{i+j} \d_i\d_{j-1}     
    &\\
      &\qquad=
    \sum_{i \leq j, j \leq p+1}
      (-1)^{i+j} \d_i\d_{j}     
    &
      &=
    \sum_{j \leq i, i \leq p+1}
      (-1)^{i+j} \d_j\d_{i}.
  \end{align*}
  Therefore,
  \begin{align*}
    \d\d 
      &=
    \sum_{0 \leq j \leq p+2} 
        (-1)^j \d_j(
                    \sum_{0 \leq i \leq p+1} (-1)^i \d_i
                   )   
      =
    \sum_{0 \leq j \leq p+2}
    \sum_{0 \leq i \leq p+1}
      (-1)^{i+j} \d_j\d_i     
    \\
      &\qquad=
    \sum_{j \leq i, i \leq p+1}
      (-1)^{i+j} \d_j\d_i
      +
    \sum_{i < j, j \leq p+2}
      (-1)^{i+j} \d_j\d_i     
      =
    0.
  \end{align*}
\end{proof}
\end{prop}

When $(G,\tht)$ is abelian, 
    Proposition \ref{prop:BMBMeq0} allows 
    the cohomology groups 
    \begin{equation*}
        H_\Gam^p(\cU,X,(G,\tht))
    \end{equation*}
    of the complex 
        $(K_\Gam^{\bul}(\cU,X,(G,\tht)),\d)$
    to be defined.
    The restriction to abelian $\Gam$-groups
    is neccesary to ensure that 
    the coboundary maps $\d^p$ 
    are group homomorphisms.
    In order to obtain a cohomology theory
    which is independent of the cover $\cU$,
    the direct limit of these cohomology 
    groups will be taken with respect
    to refinements of the cover.
    A refinement of $\cU$ consists of 
    another cover $\cV$ indexed by some set $B$,
    and a refining map 
        $r: B \rightarrow A$
    such that 
        $V_b \subset U_{r(b)}$
    for all $b \in B$.
    Such a refinement induces a refinement
    of the associated simplicial covers,
    and restriction homomorphisms 
            $
                r_*
                    :
                K^p_\Gam(\cU,X, (G,\tht)) 
                    \rightarrow      
                K^p_\Gam(\cV,X, (G,\tht))
            $
    defined by
    \begin{equation*}
        (r_*\phi)_{( b_0, \ldots, b_{p} )}
        :=
        \phi_{( r(b_0), \ldots, r(b_{p}) )}
            |_{ V_{( b_0, \ldots, b_{p} )} }.
    \end{equation*}
    These restriction homomorphisms,
    in turn, induce maps
    \begin{equation*}
        H_\Gam^p(\cU,X,(G,\tht))
            \rightarrow
        H_\Gam^p(\cV,X,(G,\tht))
    \end{equation*}
    on the cohomology of the complexes.
    In order for the direct limit 
    of cohomology groups to be well-defined,
    the maps induced on cohomology by
    two different refining maps need to be equal.
    This is true in the 
    equivariant setting,
    and in the semi-equivariant setting 
    it just needs to be checked that the 
    twisting of the coboundary map using $\tht$ 
    doesn't cause any problems.
\begin{lem}
\label{lem:SCrefine}
Let $(\cV,r)$ and $(\cV,s)$ be refinements of $\cU$
    with refining maps 
        $r,s: B \rightarrow A$.
    The maps induced on semi-equivariant cohomology
    by $r$ and $s$ are identical.
\begin{proof}
By analogy with the proof in the non-equivariant case
    (see for example
    \cite[pp.~78-79]{Raeburn.Williams.1998.Morita-equivalence-and-continuous-trace-C-star-algebras}),
    a cochain homotopy
\begin{equation*}
\xymatrix{
            & 
            K^p_\Gam(\cU,X, (G,\tht)) 
                \ar[dl]_{h^p} 
                \ar@<-0.5ex>[d]_{r_*} 
                \ar@<0.5ex>[d]^{s_*}
                \ar[r]^{\d^p} & 
            K^{p+1}_\Gam(\cU,X, (G,\tht))
                \ar[dl]_{h^{p+1}}
            \\
            K^{p-1}_\Gam(\cV,X, (G,\tht))
                \ar[r]^{\d^{p-1}} &
            K^p_\Gam(\cV,X, (G,\tht)). &
         }
\end{equation*}
    is defined by
    \begin{equation*}
        (h^p\phi)_{(b_0, \ldots, b_{p-1})}
            =
        \sum_{k=0}^{p-1}
            (-1)^k
            \phi_{(
                    r(b_0), 
                    \ldots, 
                    r(b_k), 
                    s(b_k), 
                    \ldots, 
                    s(b_{p-1})
                 )}
                \circ 
            e_k,
    \end{equation*}
    where $e_k$ is the $k$th degeneracy map.
Just as in the non-equivariant case, 
    expanding the expression 
    \begin{equation*}
        (h^{p+1} \d^p \phi)_{(b_0, \ldots, b_{p})} 
            -
        (\d^{p-1} h^{p} \phi)_{(b_0, \ldots, b_{p})} 
            \in 
        K^p_\Gam(\cV,X, (G,\tht))
    \end{equation*}
    results in a large amount of cancelation.
    The remaining expression is
    \begin{align*}
        (\d^p_0\phi)_{(
                       r(b_0), 
                       s(b_0), 
                       \ldots, 
                       s(b_p)
                     )}
            \circ 
        e_0 
            -
        (\d^p_{p+1}\phi)_{(
                       r(b_0), 
                       \ldots, 
                       r(b_p), 
                       s(b_p)
                     )}
            \circ 
        e_p. 
    \end{align*}
    The twisted coboundary maps
        $\d^0_0$ and $\d^p_{p+1}$
    involve the $\Gam$-actions 
        $\tht$ on $G$ and 
        $\sig$ on $X$,
    respectively.
    However, in the above expression, 
    the degeneracy maps $e_0$ and $e_p$ 
    ensure that $\tht$ and $\sig$
    only ever act via 
    the identity element of $\Gam$.
    Thus, the above expression simplifies to
    \begin{align*}
        \phi_{( s(b_0), \ldots, s(b_p))}
            -
        \phi_{( r(b_0), \ldots, r(b_p))}
            =
        (s_* \phi)_{( b_0, \ldots, b_p)}
            -
        (r_* \phi)_{( b_0, \ldots, b_p
             )}.
    \end{align*}
    Therefore, if 
    $\phi \in H^p_\Gam(\cV,X, (G,\tht))$
    is a cocycle,
    then
    \begin{equation*}
        (s_*\phi)
            -
        (r_*\phi)
            =
        h^{p+1} \circ \d^p(\phi)
            -
        \d^{p-1} \circ h^{p}(\phi)
            =
        \d^{p-1} \circ h^{p}(\phi),
    \end{equation*}
    which is a coboundary.
    Thus, 
    $r_*$ and $s_*$ induce the same cohomology groups.
\end{proof}   
\end{lem}

It is now possible to define 
    the semi-equivariant cohomology groups.
\begin{defn}{}
\label{def:}
The \emph{(smooth) 
    $\Gam$-semi-equivariant 
    \v{C}ech cohomology groups}
    with coefficients in 
    an abelian $\Gam$-group
    $(G,\tht)$
    are defined by
    \begin{equation*}
        H_\Gam^p(X,(G,\tht)) 
            := 
        \lim_{\rightarrow} H_\Gam^p(\cU,X,(G,\tht)),
    \end{equation*}
    where 
        $H_\Gam^p(\cU,X,(G,\tht))$
    are the cohomology groups 
    of the complex 
    $(K_\Gam^{\bul}(\cU,X,(G,\tht)),\d)$, and 
    the direct limit is taken 
    with respect to refinements of $\cU$.
\end{defn}

Semi-equivariant cohomology is functorial 
    with respect to homomorphisms 
    of abelian $\Gam$-groups.
\begin{prop}
\label{prop:Gm-gives-Cm}
A homomorphism $\al: A \rightarrow B$ 
    of abelian $\Gam$-groups 
    induces a morphism of complexes
    \begin{equation*}
        \al^\bul: (K_\Gam^\bul(\cU,X,A),\d) \rightarrow (K_\Gam^\bul(\cU,X,B),\d)
    \end{equation*}
    defined by $(\al^p\phi)_{a^p} := \al \circ \phi_{a^p}.$
\begin{proof}
    Let $\tht$ be the $\Gam$-action on $A$ and $\vtht$ be the $\Gam$-action on $B$.
    As $\al$ is a homomorphism of $\Gam$-groups
    $\al^p \circ \tht^{x^p}_i = \vtht^{x^p}_i \circ \al^p $
    for all $x^p \in X^p$ and $0 \leq i \leq p$.
    Thus,
    \begin{align*}
    (\al^{p+1} (\d_i\phi))_{a^{p+1}}(x^{p+1})
      &=
    \al \circ (\d_i\phi)_{a^{p+1}}(x^{p+1}) \\
      &=
    \al \circ 
      \tht^{x^{p+1}}_i \circ 
      \phi_{d_i(a^{p+1})} \circ 
      d_i(x^{p+1}) \\
      &=
    \vtht^{x^{p+1}}_i \circ 
      \al \circ 
      \phi_{d_i(a^{p+1})} \circ 
      d_i(x^{p+1}) \\
      &=
    \vtht^{x^{p+1}}_i \circ 
      (\al^p\phi)_{d_i(a^{p+1})} \circ 
      d_i(x^{p+1}) \\
      &=
    (\d_i (\al^p \phi))_{a^{p+1}}(x^{p+1}).
    \end{align*}
    Therefore,  
    $\al^{p+1} \circ \d = \d \circ \al^p$ 
    and $\al^p$ defines a morphism of complexes.
\end{proof}
\end{prop}

Given a short exact sequence of abelian $\Gam$-groups,
    connecting maps for a long exact sequence can be constructed.
\begin{thm}\label{thm:LESab}
A short exact sequence of abelian $\Gam$-groups
    \begin{equation*}
        1 \rightarrow   
        A \os{\al}{\rightarrow}
        B \os{\be}{\rightarrow}
        C \rightarrow
        1
    \end{equation*}
  induces a long exact sequence
  \begin{equation*}
    \ldots
    \os{\Del^{p-1}}{\rightarrow}
  H_\Gam^p(X,A)     
    \os{\al^p}{\rightarrow}
  H_\Gam^p(X,B) 
    \os{\be^p}{\rightarrow} 
  H_\Gam^p(X,C)
    \os{\Del^{p}}{\rightarrow}
  H_\Gam^{p+1}(X,A)
    \os{\al^{p+1}}{\rightarrow}
    \ldots,
  \end{equation*}
  where
        $\Del^p(\phi) := [\d(\psi)]$
  for any element 
    $\psi \in K_\Gam^p(B)$
  such that
    $\be^p(\psi) = \phi$.
\begin{proof}
    The proposition follows by standard diagram chasing arguments
    applied to the exact sequence of complexes
    \begin{equation*}
            1 
                \rightarrow 
            (K_\Gam^\bul(X,A),\d)
                \os{\al^\bul}{\rightarrow} 
            (K_\Gam^\bul(X,B),\d)
                \os{\be^\bul}{\rightarrow} 
            (K_\Gam^\bul(X,C),\d)
                \rightarrow
            1.
    \end{equation*}
  For an example, see the proof of 
  \cite[Theorem 4.30]{Raeburn.Williams.1998.Morita-equivalence-and-continuous-trace-C-star-algebras}.
\end{proof}
\end{thm}

\section{The Semi-equivariant Dixmier-Douady Invariant}
\label{sec:semiDix}

In order to apply semi-equivariant cohomology 
    to 
    the classification 
    of semi-equivariant liftings,
    its relationship with 
    semi-equivariant principal bundles 
    must be clarified.
    By Theorem \ref{thm:PB-TC-bij}, 
    this reduces 
    to the problem of 
    relating semi-equivariant transition cocycles
    and
    semi-equivariant cohomology classes.
In this section, 
    semi-equivariant transition cocycles 
    will be interpreted 
    as degree-$1$ cocycles
    which can take values in a 
    non-abelian coefficient group.
    An analogue of Theorem \ref{thm:LESab} 
    will be proved 
    that constructs a connecting map 
    from the transition cocycles 
    into degree-$2$ cohomology.
    The theorem can be used 
    to classify certain liftings 
    of semi-equivariant principal bundles
    between non-abelian structure groups.
    This method has its origins in the
    work of Dixmier-Douady 
    on continuous trace $C^*$-algebras
\cite{Dixmier-Douady.1963.Champs-continus-despaces-Hilbertiens-et-de-C-star-algebres}.
    See also 
    \cite[\S4]{2008.Brylinski.Loop-spaces-characteristic-classes-and-geometric-quantization}
    and
    \cite[\S4.3]{Raeburn.Williams.1998.Morita-equivalence-and-continuous-trace-C-star-algebras}.

To begin, 
    note that the $p$-cochains 
    of Definition \ref{def:cochains} and 
    the twisted pullback maps of Definition \ref{def:TPB}
    are well-defined for non-abelian $\Gam$-groups.
    Thus, 
    it is possible to make the following definitions.
\begin{defn}\label{def:TCnabH}
    \begin{align}
        \TCocycles_\Gam^0(\cU,X&,(G,\tht))
            := 
        \set{
                \mu \in  K_\Gam^0(\cU,X,(G,\tht))
                \mid
                (\d_1 \mu)^\inv (\d_0 \mu) = 1
            } \label{def:TCnabH0}
        \\
        \TCocycles_\Gam^1(\cU,X&,(G,\tht))
        \nonumber
        \\
            &:= 
        \set{
                \phi \in  K_\Gam^1(\cU,X,(G,\tht))
                \mid
                (\d_1\phi)^\inv (\d_2\phi) (\d_0\phi) = 1
            }  / \sim \label{def:TCnabH1}   
    \end{align} 
    where $\phi^1 \sim \phi^2$ 
    if and only if 
    there exists a 
    $\mu \in  K_\Gam^0(\cU,X,(G,\tht))$
    such that
        $(\d_1\mu) \phi^1 = \phi^2 (\d_0\mu)$.
\end{defn}

The set 
    $\TCocycles_\Gam^1(\cU,X,(G,\tht))$ 
    is just 
    $\TCocycles_\Gam^\iso(\cU,X,(G,\tht))$
    with 
    the transition cocycle condition and 
    equivalence condition
    expressed in terms of twisted pullback maps.
    Note that the particular order of the terms 
        $\d_i\mu$ in \eqref{def:TCnabH0} and 
        $\d_i\phi$ in \eqref{def:TCnabH1}
    is important 
    as the elements $\mu$ and $\phi$ take values in $G$,
    which is not necessarily abelian.
    When $G$ is abelian, 
    these terms may be rearranged 
    to give the corresponding cocycle properties 
    in semi-equivariant cohomology.
    An abelian structure group also ensures that 
    pointwise multiplication 
    is a well-defined group structure on
        $\TCocycles_\Gam^0$ and $\TCocycles_\Gam^1$,
    which, in general, are only pointed sets.
\begin{thm}\label{thm:H1isTC}
When $G$ is abelian
    \begin{align}
        \TCocycles_\Gam^0(\cU,X,(G,\tht)) 
            &\iso 
        H_\Gam^0(\cU,X,(G,\tht))  \label{eq:TC0H0}
    \\
        \TCocycles_\Gam^1(\cU,X,(G,\tht))
            &\iso 
        H_\Gam^1(\cU,X,(G,\tht)). \label{eq:TC1H1}
    \end{align}
\begin{proof}
When $G$ is abelian, 
    the defining condition on
        $\TCocycles_\Gam^0(\cU,X,(G,\tht))$
    and the $0$-cocycle condition 
    on cohomology are equivalent as
    \begin{equation*}
        0 = -(\d_1 \mu) +  (\d_0 \mu) 
          = (\d_0 \mu) - (\d_1 \mu)
          = \d \mu.
    \end{equation*}
    This proves \eqref{eq:TC0H0}.
Similarly, the defining condition on
        $\TCocycles_\Gam^1(\cU,X,(G,\tht))$ 
    and the $1$-cocycle condition on cohomology 
    are equivalent as
    \begin{equation*}
        0 = -(\d_1\phi) + (\d_2\phi) + (\d_0\phi) = (\d_0\phi) - (\d_1\phi) + (\d_2\phi)  = \d\phi,
    \end{equation*}
    and the equivalence relations on 
        $\TCocycles_\Gam^1(\cU,X,(G,\tht))$ 
    and
        $H_\Gam^0(\cU,X,(G,\tht))$ 
    are the same as
    \begin{align*}
        (\d_1\mu) + \phi^1 &= \phi^2 + (\d_0\mu) 
        \\
         \phi^1 - \phi^2   &=  (\d_0\mu) - (\d_1\mu) 
        \\
         \phi^1 - \phi^2   &=  \d\mu.
    \end{align*}
    These two facts imply \eqref{eq:TC1H1}.
\end{proof}
\end{thm}

Together, 
        Theorem \ref{thm:LESab} and 
        Theorem \ref{thm:H1isTC} 
    enable liftings 
    of semi-equivariant principal bundles
    between abelian structure groups 
    to be classified. 
    However, 
    the construction of a Dirac operator 
    involves the construction of liftings
    between non-abelian groups.
The next theorem is a generalisation 
    of Theorem \ref{thm:LESab}
    that can be used 
    to classify certain liftings 
    between non-abelian structure groups.
\vspace{4em}
\begin{thm}
\label{thm:LESnab}
A short exact sequence of $\Gam$-groups
    \begin{equation*}
        1 \rightarrow   
        A \os{\al}{\rightarrow}
        B \os{\be}{\rightarrow}
        C \rightarrow
        1,
    \end{equation*}
    where $\al(A)$ is central in $B$,
    induces an exact sequence of pointed sets
    \begin{align*}
      0
      &\rightarrow
    H_\Gam^0(X,A)     
      \os{\al^0}{\rightarrow}
    \TCocycles_\Gam^0(X,B)
      \os{\be^0}{\rightarrow}
    \TCocycles_\Gam^0(X,C)
    \os{\Del^{0}}{\rightarrow}
      \ldots
    \\
    \ldots& 
    \os{\Del^{0}}{\rightarrow}
    H_\Gam^1(X,A)     
      \os{\al^1}{\rightarrow}
    \TCocycles_\Gam^1(X,B)
      \os{\be^1}{\rightarrow}
    \TCocycles_\Gam^1(X,C)
      \os{\Del^{1}}{\rightarrow}
    H_\Gam^2(X,A),
    \end{align*}
    where
    the connecting maps $\Del^0$ and $\Del^1$
    are defined by
    \begin{align*}
    \Del^0([\mu]) 
        &:= [(\d_1\eta)^\inv (\d_0\eta)]
    &
    \Del^1([\phi]) 
        := 
    [(\d_1\psi)^\inv (\d_2\psi) (\d_0\psi)],
    \end{align*}
    for any
    $\eta \in K_\Gam^0(X,B)$,
    $\psi \in K_\Gam^1(X,B)$
    satisfying
    $\be^0(\eta) = \mu$,
    $\be^1(\psi) = \phi$.
\begin{proof}
The diagram chasing arguments used in the proof of Theorem \ref{thm:LESab} do not apply directly.
    However, they can be imitated 
    while carefully working around 
    any lack of commutivity in the groups $B$ and $C$.
Note that 
    Proposition \ref{prop:Gm-gives-Cm} and 
    Proposition \ref{prop:tpb-si} continue to hold 
    when the structure groups involved are non-abelian.
    Thus, the twisted pullback maps $\d_i$ 
    commute with 
    the maps $\al^i$ and $\be^i$ 
    induced by $\al$ and $\be$,
    and also satisfy the simplicial identity 
    $\d_j \circ \d_i = \d_i \circ \d_{j-1}$ for $i < j$.

First, the map $\Del^0$ will be considered.
    Let 
        $
        \nu 
            := 
        (\d_1\eta)^\inv (\d_0\eta) \in K^1_\Gam(X,B)
        $. 
    The cochain $\eta$ is a lifting by $\be$ of $\mu$ so
        $\be(\nu) = 1$.
    Thus, $\nu$ takes values in $\ker(\be) \iso A$ and defines an element of $K_\Gam^1(X,A)$.
    The simplicial identity can be used to show that the cochain $\nu$ satisfies the cocycle property,
    \begin{align*}
       (\d_1\nu)^\inv (\d_0\nu) 
            &=
        \Big[ (\d_1\d_1\nu)^\inv (\d_1\d_0\nu) \Big]^\inv
        \Big[ (\d_0\d_1\nu)^\inv (\d_0\d_0\nu) \Big] \\
            &=
         (\d_1\d_0\nu)^\inv (\d_1\d_1\nu)
         (\d_0\d_1\nu)^\inv (\d_0\d_0\nu) \\
            &=
         (\d_1\d_0\nu)^\inv (\d_1\d_1\nu)
         (\d_0\d_1\nu)^\inv (\d_0\d_0\nu)  \\
            &=
         (\d_0\d_0\nu)^\inv (\d_1\d_1\nu)
         (\d_1\d_1\nu)^\inv (\d_0\d_0\nu)  \\
            &=
        1.
    \end{align*}
    Therefore, 
    $\Del^0([\mu]) := [\nu] \in H_\Gam^1(X,A)$.
Next, it needs to be shown that 
        $\Del^0([\mu]) := [(\d_1\eta)^\inv (\d_0\eta)]$
    is independent of the choice of $\eta$.
    Let $\eta' \in K^0_\Gam(X,B)$ be another element 
    such that $\be(\eta') = \mu$. 
    Set $\omg := \eta'\eta^\inv$ and 
        $\nu' := (\d_1\eta')^\inv (\d_0\eta') 
              \in K^1_\Gam(X,B)$.
    Then $\be(\omg) = \be(\eta'\eta^\inv) 
                    = \mu\mu^\inv
                    = 1$.
    Thus, $\omg$ defines an element of $K_\Gam^0(X,A)$
    and $\d\omg \in K_\Gam^1(X,A)$ is a coboundary.
    Using the fact that $\nu$ and $\d\omg$ take values in the abelian group $A$,
    \begin{align*}
        (\d\omg)\nu
            &=
        (\d\omg)(\d_1\eta)^\inv(\d_0\eta) \\
            &=
        (\d_1\eta)^\inv(\d\omg)(\d_0\eta) \\
            &=
        (\d_1\eta)^\inv 
        (\d_1\eta)
        (\d_1\eta')^\inv
        (\d_0\eta')
        (\d_0\eta)^\inv
        (\d_0\eta) \\
            &=
        (\d_1\eta')^\inv
        (\d_0\eta') \\
            &=
        \nu'.
    \end{align*}
    Therefore, 
    $[\nu] = [\nu'] \in H_\Gam^1(X,A)$.

In order to examine the map $\Del^1$, 
    let $\nu := (\d_1\psi)^\inv (\d_2\psi) (\d_0\psi) 
             \in K^2_\Gam(X,B)$.
    The cochain $\psi \in K_\Gam^1(X,B)$ 
    is a $\be$-lifting 
    of the cocycle 
        $\phi \in \TCocycles_\Gam^1(X,C)$ 
    so 
        $\be(\nu) = 1$.
    Therefore, $\nu$ defines an element of $K_\Gam^2(X,A)$.
Using the simplicial identity,
    and the fact that 
    $\nu$ takes values in the centre of $B$,
    it can be shown that 
    $\nu$ satisfies the $2$-cocycle propery.
    First, compute
    \begin{align*}
       (\d_1 \nu)(\d_3 \nu) 
        &= (\d_1\d_1 \psi)^\inv (\d_1\d_2 \psi)  (\d_1\d_0 \psi)        (\d_3 \nu)                                              \\ 
        &= (\d_1\d_1 \psi)^\inv (\d_1\d_2 \psi)  (\d_3 \nu)             (\d_1\d_0 \psi)                                         \\ 
        &= (\d_1\d_1 \psi)^\inv (\d_1\d_2 \psi)  \Big[ (\d_3\d_1 \psi)^\inv   (\d_3\d_2 \psi)         (\d_3\d_0 \psi) \Big]   (\d_1\d_0 \psi)\\ 
        &= (\d_1\d_1 \psi)^\inv (\d_1\d_2 \psi)  \Big[ (\d_1\d_2 \psi)^\inv   (\d_3\d_2 \psi)         (\d_3\d_0 \psi)  \Big] (\d_1\d_0 \psi)\\ 
        &= (\d_1\d_1 \psi)^\inv (\d_3\d_2 \psi)  (\d_3\d_0 \psi)        (\d_1\d_0 \psi)                                         \\ 
        &= (\d_1\d_1 \psi)^\inv (\d_3\d_2 \psi)  \Big[ (\d_2\d_0 \psi)  (\d_2\d_0 \psi)^\inv \Big] (\d_3\d_0 \psi)  (\d_1\d_0 \psi) \\ 
        &= (\d_2\d_1 \psi)^\inv (\d_2\d_2 \psi)  \Big[ (\d_2\d_0 \psi)  (\d_0\d_1 \psi)^\inv \Big] (\d_0\d_2 \psi)  (\d_0\d_0 \psi) \\ 
        &= \Big[ (\d_2\d_1 \psi)^\inv (\d_2\d_2 \psi)   (\d_2\d_0 \psi) \Big] \Big[ (\d_0\d_1 \psi)^\inv  (\d_0\d_2 \psi)  (\d_0\d_0 \psi) \Big] \\ 
        &= (\d_2 \nu)(\d_0 \nu).
    \end{align*}
    Then
    \begin{align*}
        (\d\nu) 
            &= (\d_0 \nu)(\d_1 \nu)^\inv(\d_2 \nu)(\d_3 \nu)^\inv \\
            &= (\d_0 \nu)(\d_2 \nu)(\d_3 \nu)^\inv(\d_1 \nu)^\inv \\
            &= (\d_0 \nu)(\d_2 \nu) \Big[ (\d_1 \nu)(\d_3 \nu) \Big]^\inv \\
            &= (\d_0 \nu)(\d_2 \nu) \Big[ (\d_0 \nu)(\d_2 \nu) \Big]^\inv \\
            &= 1,
    \end{align*}
    and so $[\nu] \in H_\Gam^2(X,A)$.

Next, 
    it needs to be shown that $\Del^1$ is well-defined. 
    Specifically, that
    \begin{equation*}
        \Del^1([\phi]) 
            := 
        [(\d_1\psi)^\inv (\d_2\psi) (\d_0\psi)]
    \end{equation*}
        is independent of the choice of $\psi$, 
    and
        depends only on 
        the class of 
        $\phi$ in $\TCocycles_\Gam^1(X,C)$.
To prove the first statement, let 
        $\psi' \in K_\Gam^1(X,B)$ 
    be another $\be$-lifting of $\phi$ and 
        $\nu' := (\d_1\psi')^\inv (\d_2\psi') (\d_0\psi')$ 
    be the corresponding element of $H^2_\Gam(X,A)$.
    If $\omg := \psi'\psi^\inv$ then 
        $
            \be(\omg) 
                = 
            \be(\psi'\psi^\inv) 
                = 
            \phi\phi^\inv 
                = 
            1.
        $
    Thus, $\omg \in K_\Gam^1(X,A)$ and $\d\omg \in K^2_\Gam(X,A)$ is a coboundary .
    Next, using the fact that $\omg$ takes values in the center of $B$,
    \begin{align*}
    (\d\omg)\nu 
        &=  
    (\d_0\omg)(\d_1\omg)^\inv (\d_2\omg)              
    (\d_1\psi)^\inv (\d_2\psi) (\d_0\psi)\\
        &=  
    (\d_1\psi)^\inv(\d_1\omg)^\inv                  
    (\d_2\omg)(\d_2\psi)                      
    (\d_0\omg)(\d_0\psi)\\
        &=  
    (\d_1\psi)^\inv(\d_1\psi'\psi^\inv)^\inv       
    (\d_2\psi'\psi^\inv)(\d_2\psi)           
    (\d_0\psi'\psi^\inv)(\d_0\psi)\\
        &=  
    (\d_1\psi)^\inv(\d_1\psi) (\d_1\psi')^\inv     
    (\d_2\psi') (\d_2\psi)^\inv(\d_2\psi)    
    (\d_0\psi')(\d_0\psi)^\inv(\d_0\psi)\\
        &=  
    (\d_1\psi')^\inv (\d_2\psi') (\d_0\psi') \\
        &=  
    \nu'.
    \end{align*}
    Therefore, $[\nu] = [\nu'] \in H_\Gam^2(X,A)$.

In order to prove that $\Del^1([\phi])$ 
    depends only on the class of $\phi$,
    suppose that $\phi$ is a coboundary i.e. that $\phi = (\d_1\til{\phi})^\inv (\d_0\til{\phi})$ for some $\til{\phi} \in K_\Gam^0(X,C)$.
    By surjectivity of $\be$, there exists an element $\til{\psi}$ such that $\be(\til{\psi})= \til{\phi}$. 
    Then $\psi := (\d_1\til{\psi})^\inv (\d_0\til{\psi})$ is a lifting by $\be$ of $\phi$ as
    \begin{align*}
        \be(\psi) 
            &=  \be\Big[ (\d_1\til{\psi})^\inv (\d_0\til{\psi}) \Big]  \\
            &=  (\be\d_1\til{\psi})^\inv (\be\d_0\til{\psi}) \\
            &=  (\d_1\be\til{\psi})^\inv (\d_0\be\til{\psi}) \\
            &=  (\d_1\til{\phi})^\inv (\d_0\til{\phi}) \\
            &=  \phi. 
    \end{align*}
    So, again applying the simplicial identity,
    \begin{align*}
        \Del^1([\phi]) 
            &= [(\d_1\psi)^\inv (\d_2\psi) (\d_0\psi)]  \\
            &= [ (\d_1\d_0\til{\psi})^\inv  (\d_1\d_1\til{\psi})  (\d_2\d_1\til{\psi})^\inv (\d_2\d_0\til{\psi})  (\d_0\d_1\til{\psi})^\inv (\d_0\d_0\til{\psi})] \\
            &= [ (\d_0\d_0\til{\psi})^\inv  (\d_1\d_1\til{\psi})  (\d_1\d_1\til{\psi})^\inv (\d_0\d_1\til{\psi})  (\d_0\d_1\til{\psi})^\inv (\d_0\d_0\til{\psi})] \\
            &= 1.
    \end{align*}
    Thus, $\Del^1([\phi])$ depends only on the class of $\phi$ in $\TCocycles_\Gam^1(X,C)$.
\end{proof}
\end{thm}

It is now possible to 
    define the
    semi-equivariant Dixmier-Douady invariant
    and resolve the main problem 
    of this paper.
\begin{defn}
\label{defn:SEDD}
The 
    \emph{semi-equivariant Dixmier-Douady invariant} 
    of a 
    $\Gam$-semi-equivariant 
    principal $C$-bundle $P$
    associated to a central exact sequence
    \begin{equation*}
        1 \rightarrow   
        A \os{\al}{\rightarrow}
        B \os{\be}{\rightarrow}
        C \rightarrow
        1
    \end{equation*}
    is defined by
    \begin{equation*}
    DD(P) := \Del^1([\phi]) \in H^2_\Gam(X,A),
    \end{equation*}
    where 
    $\Del^1$ is the connecting map 
    provided by Theorem \ref{thm:LESnab} and
    $[\phi]$ is the transition cocycle 
    associated to $P$ by Proposition \ref{prop:TC-from-PB}.
\end{defn}

\begin{thm}
\label{thm:mainthm}
The exact sequence produced 
    by Theorem \ref{thm:LESnab}
    implies that
\begin{enumerate}
\item 
    $P$ can be lifted
    by $\be$ 
    if and only if
    $DD(P) = 0$,
\item 
    when $DD(P)=0$,
    the liftings of $P$ by $\be$ 
    correspond non-canonically
    to the classes of
    $H_\Gam^1(X,A)$.
\end{enumerate}
\end{thm}

\section{Related Work and Applications}
\label{sec:relwork}

Semi-equivariant \v{C}ech cohomology
    $H_\Gam^\bul(X,(G,\tht))$
    is closely related 
    to several other cohomology theories.
    For example,
\begin{enumerate}
\item 
    When $\Gam$ is the trivial group,
    $H_\Gam^\bul(X,(G,\tht))$
    is \v{C}ech cohomology.
\item 
    When $\tht$ is the trivial action,
    $H_\Gam^\bul(X,(G,\tht))$
    is equivariant \v{C}ech cohomology
    $\check{H}_\Gam^\bul(X,G)$.
    When $X$ is a compact manifold
    acted upon by a finite group,
    the equivariant \v{C}ech cohomology
    can be related to
    Grothendieck's equivariant sheaf cohomology
    \cite[\S 5.5]{Grothendieck.1957.Sur-quelques-points-d-algebre-homologique}
    or
    Borel cohomology
    \cite[\S A]{Brylinski.2000.Gerbes-on-complex-reductive-Lie-groups},
    \cite[\S 3.3]{Gomi.2005.Equivariant-smooth-Deligne-cohomology}.
    
    Note that 
    there is a restriction
    homomorphism
    \begin{equation*}
    H_\Gam^p(X,(G,\tht))
    \rightarrow
    H_{\Gam_G}^p(X,(G,\tht))
    \iso
    \check{H}_{\Gam_G}^p(X,G),
    \end{equation*}
    where
    $\Gam_G \subseteq \Gam$
    is the stabiliser subgroup
    that acts trivially on $G$.
    In this way,
    the semi-equivariant cohomology
    can be regarded
    as a restriction of 
    equivariant cohomology.
\item 
    When $X$ is a point,
    $H_\Gam^\bul(X,(G,\tht))$
    is the group cohomology
    $H^\bul(\Gam,G_\tht)$
    of $\Gam$
    with coefficients in the $\Gam$-module
    $G_\tht$ defined by $G$ and $\tht$
    \cite[p.~35]{Berhuy.2010.An-introduction-to-Galois-cohomology-and-its-applications}.
    With this in mind,
    semi-equivariant cohomology
    can be viewed as a 
    cross between 
    group cohomology and 
    equivariant cohomology.

\item
    When $X$ is a Real space and 
    $\kap$ is the conjugation action 
    on $\Unitary(1)$,
    $H_{\Gal(\C/\R)}^\bul(X,(\Unitary(1),\kap))$
    is closely related 
    to Real \v{C}ech cohomology of
    \cite{Moutuou.2013.On-groupoids-with-involution-and-their-cohomology},
    and the Real sheaf cohomology
    defined in 
    \cite{Hekmati-Murray-Szabo-Vozzo.2016.Real-bundle-gerbes-orientifolds-and-twisted-KR-homology}.
    Note that, in this case,
    the semi-equivariant cohomology
    incorporates aspects of 
    equivariant \v{C}ech cohomology
    and
    Galois cohomology for the 
    field extension $\C/\R$.
\end{enumerate}

An important application 
    of Theorem \ref{thm:LESnab}
    arises in the study of $\Spinc$-structures
    on Real spaces 
    \cite{Atiyah.1966.K-Theory-and-Reality}
    and 
    orientifolds
    \cite{Freed-Moore.2013.Twisted-equivariant-matter}.
    This is the original motivation 
    for the present paper.
    Such structures correspond to 
    semi-equivariant liftings of 
    equivariant principal 
    $\SOrth(n)$-bundles
    via the central exact sequence
    \begin{equation*}
        1 
            \rightarrow  
        (\Unitary(1),\kap \circ \eps) 
            \rightarrow 
        (\Spinc(n),\kap \circ \eps) 
            \os{\Ad^c}{\rightarrow} 
        (\SOrth(n),\id \circ \eps)
            \rightarrow 
        1.
    \end{equation*}
    Here $\eps: \Gam \rightarrow \Z_2$
    is a homomorphism from a finite group $\Gam$,
    and
    $\kap$ denotes the conjugation action
    on $\Spinc(n)$ and $\Unitary(1)$.
    The topic of $\Spinc$-structures 
    on orientifolds 
    will be treated in a forthcoming paper.



\bibliographystyle{plain}
\bibliography{/home/simon/Library/Library-Contents}


\end{document}